\documentclass[11.5pt]{amsart}
\usepackage{amsfonts}
\usepackage{amssymb}
\usepackage{pict2e}
\usepackage[all,poly,graph]{xy}

\font\got=eufm10
\font\gots=eufm10 at 7pt

\def\g2{ \hbox{\got g}_2}
\def\f4{\hbox{\got f}_4}
\def\d4{\hbox{\got d}_4}
\def\e6c{\hbox{\got e}_6^\mathbb{C}}
\def\fs4{\hbox{\gots f}_4}
\def\F4{\hbox{\got F}_4}
\def\Fs4{\hbox{\gots F}_4}

\def\es7{ \hbox{\got e}_7}
\def\es8{ \hbox{\got e}_8}
\def\id{\mathop{\hbox{\rm id}}}

\def\Ok{{\mathcal O} }

\def\F{\mathbb F}
\def\R{\mathbb R}
\def\imag{\textbf{i}}
\newcommand{\ii}{\textbf{i}}
\def\aut{\mathop{\rm Aut}}

\def\der{\mathop{\rm Der}}
\def\End{\mathop{\rm End}}

\def\ad{\mathop{\rm ad}}

\def\Mat{\mathop{\rm Mat}}

\def\span#1{\langle #1 \rangle}

\def\tr{\mathop{\rm tr}}

 \def\e6{  \mathfrak{e}_6 }
 \def\sig#1{ \mathfrak{e}_{6,#1} }

\title[Satake diagrams on $\e6$  ]{ On the real forms of the exceptional Lie algebra $\e6$ and their Satake diagrams}  %%%para que no salga
\author[C. Draper]{Cristina Draper}
\thanks{\ Partially supported by MCYT grant   MTM2013-41208-P  and by   the Junta de Andaluc\'{\i}a PAI
projects  FQM-336,  FQM-2467, FQM-3737.}
\address{Cristina Draper Fontanals: Departamento de
Matem\'atica Aplicada\\ Escuela de las Ingenier\'{\i}as\\ Ampliaci\'{o}n Campus de Teatinos, S/N, 29071 M\'alaga,
Spain\\cdf@uma.es}
\author[V. Guido]{Valerio Guido}
\thanks{\  Partially supported by   the Ministero dell'Istruzione, dell'Universit\`{a} e della
Ricerca (MIUR) with a grant for the Ph.D. at the Universit\`{a} del Salento}
\address{Valerio Guido: Dpto. di Matematica E. De Giorgi\\ Universit\`{a} del Salento\\
Via Provinciale Lecce-Arnesano, 73100, Lecce, Italy\\valerio.guido@unisalento.it}

\usepackage{graphicx}
\begin{document}

\setlength{\unitlength}{0.06in}

\maketitle

%%%%%%%%%%%%%%%%%%%%%%%%%%%%%%%%%%%%%%%%%%%%%%%%%%%%%%%%%%%%%%%%%%%%%%%%%%%%%%%%%%%%%%%%%%%%%%%%%%%%%%%%%%%%%%%%%%%%%%%%%%%%%%%%%%%%%%%

\begin{abstract}
Satake diagrams of the real forms $\sig{-26}$, $\sig{-14}$ and $\sig{2}$ are carefully developed. The first real form is constructed with an Albert algebra  and the other ones by using the two paraoctonion algebras and certain symmetric construction of the magic square.
\end{abstract}

%%%%%%%%%%%%%%%%%%%%%%%%%%%%%%%%%%%%%%%%%%%%%%%%%%%%%%%%%%%%%%%%%%%%%%%%%%%%%%%%%%%%%%%%%%%%%%%%%%%%%%%%%%%%%%%%%%%%%%%%%%%%%%%%%%%%%%%

\section{Introduction}

The real simple Lie algebras were classified by Cartan in 1914 in \cite{Cartanclasificareales}. This paper required a great amount of computations, and the classification was realized without the Cartan involution. Cartan came back to this classification in \cite{Cartan1927a}, where he identified the maximal compactly imbedded subalgebra in each case. The numbering appeared in that work is the used here in Table I. Soon he completed the classification by relating Lie algebras and geometry in \cite{CartanSatake}. This is the paper containing more information about the exceptional Lie algebras.
Several authors along the XX century provided different classifications trying to simplify Cartan's arguments. Most of them used a maximally compact Cartan subalgebra, but not all,
 Araki's approach \cite{araki} was based on choosing a maximally noncompact Cartan subalgebra. This method for classifying was a considerable improvement in clarity and simplicity.
 The classification was stated in terms of certain diagrams, called \emph{Satake  diagrams}, described by Helgason in \cite[p.\,534]{Helgasson} based on the facts developed by Satake in \cite{Satake}.
  Although Satake diagrams were contained in \cite{araki} and the restricted root systems and their multiplicities were given by  \cite{CartanSatake},
  according to Helgason's words  \cite[p.\,534]{Helgasson} Cartan stated the results for exceptional algebras   without proof.
 Since then, many textbooks contain   Satake diagrams (more historical notes can be found in \cite{knapp} and some   examples in \cite{enci,parabolic}), but, as far as we know, it is very difficult to find  details of how these diagrams were obtained in the nonclassical cases. Our objective here is to construct the diagrams starting with concrete models, which besides allows to obtain  a lot of valuable information.  %as much information as possible.

In our work in progress about gradings on the five real forms of the exceptional Lie algebra of type $\e6$, which tries to dive in the structure of such interesting simple real Lie algebras, Satake diagrams have been very useful because they codify some aspects of the structure of the corresponding real semisimple Lie algebras, namely, some  questions related to gradings   are encoded in the Satake diagram. A result in this line is the following
 \cite[Theorem~3]{Cheng}:
A simple Lie algebra $\mathfrak{g}$ admits a $\mathbb{Z}$-grading of the second kind (that is, $\mathfrak{g}=\mathfrak{g}_{-2}\oplus \mathfrak{g}_{-1}\oplus\mathfrak{g}_{0}\oplus\mathfrak{g}_{1}\oplus\mathfrak{g}_{2}  $) if and only if there is $\alpha$ a long root of $\Delta$ such that the multiplicity
$m_{\bar \alpha}=1$. But, when it was applied in \cite{tesis}   to obtain a fine grading  by the group $\mathbb{Z}^2\times \mathbb{Z}_2^3$ on $\sig{-14}$ (and on $\sig{2}$ and $\sig{6}$), it was necessary a more precise knowledge of the restricted roots than that one summarized in Table I.

The structure of this paper is as follows.
After recalling some basic facts about real forms, we explain how the Satake diagrams are constructed starting with a Cartan decomposition and a Cartan subalgebra adapted in some sense to this decomposition.
Before proceeding to compute the Satake diagrams of   $\sig{-26}$,  $\sig{-14}$  and $\sig{2}$ in Sections~\ref{sec26}, \ref{sec14} and \ref{sec2} respectively, we have enclosed a section about composition and symmetric composition algebras, Jordan algebras and constructions of exceptional Lie algebras based on these related structures, just because  the constructions we  have taken of our three real forms   make use of these nonassociative algebras.

%%%%%%%%%%%%%%%%%%%%%%%%%%%%%%%%%%%%%%%%%%%%%%%%%%%%%%%%%%%%%%%%%%%%%%%%%%%%%%%%%%%%%%%%%%%%%%%%%%%%%%%%%%%%%%%%%%%%%%%%%%%%%%%%%%%%%%%
\section{Preliminaries about real forms }

Given a real vector space $V$, we call the \emph{complexification} of $V$, and we denote it by $V^{\mathbb{C}}$,
to the complex vector space $V\otimes_\mathbb{R}\mathbb{C}=V\oplus \imag V$ ($\imag\in \mathbb{C}$ the imaginary unit).
If $ \mathfrak{g}$ is a real Lie algebra, then the complexification
$ \mathfrak{g}^\mathbb{C}$ is a complex Lie algebra with the usual extension of the bracket, that is,
$$
[x_1+\imag y_1,x_2+\imag y_2]=[x_1,x_2]-[y_1,y_2]+\imag([x_1,y_2]+[y_1,x_2])
$$
for any $x_1,x_2,y_1,y_2\in\mathfrak{g}$. In this case, $\sigma\colon\mathfrak{g}^{\mathbb{C}}\to\mathfrak{g}^{\mathbb{C}}$
given by $\sigma(x+\imag y)=x-\imag y$ is an order two antiautomorphism of $\mathfrak{g}^{\mathbb{C}}$, called the \emph{conjugation}
related to $\mathfrak{g}$. Note that $\mathfrak{g}$ coincides with the set of elements of  $\mathfrak{g}^{\mathbb{C}}$ fixed by this conjugation.

If $L$ is a complex Lie algebra and $\mathfrak{g}\subset L$ is a real  subalgebra, it is said that $\mathfrak{g}$ is a \emph{real form} of $L$ when
 $\mathfrak{g}^{\mathbb{C}}=L$. Two real forms $\mathfrak{g}_1$ and  $\mathfrak{g}_2$ of the same complex Lie algebra $L$ (with related conjugations $\sigma_1$ and $\sigma_2$ respectively) are isomorphic if and only if there is $f\in\aut(\mathfrak{g})$ such that $\sigma_2=f\sigma_1f^{-1}$.

Given $\mathfrak{g}$ a real semisimple Lie algebra,
  $\mathfrak{g}$  is said to be \emph{split} if it contains a Cartan subalgebra  $\mathfrak{h}$ such that $\ad h$ is diagonalizable
 over $\mathbb{R}$ for any $h\in\mathfrak{h}$;
and $\mathfrak{g}$ is said to be \emph{compact} if its Killing form is definite (necessarily negative). A well-known result states that any complex semisimple Lie algebra contains both a split and a compact
real form.

The importance of the real forms is due to the following. If $\mathfrak{g}$ is a simple real Lie algebra, either  $\mathfrak{g}$ is   just a complex simple Lie algebra, but considered as a
real Lie algebra,  or $\mathfrak{g}^\mathbb{C}$  is  simple, so that $\mathfrak{g}$ is  a real form of $\mathfrak{g}^\mathbb{C}$.

 The real forms of a complex simple Lie algebra $L$ are characterized by the signature of their Killing form. 
  (By abuse of notation, sometimes we speak about the signature of $L$ to refer the signature of the Killing form of $L$.)
 For instance, the signature of the split real form coincides with  the rank of $L$ and the signature of the compact real form is equal to $-\dim L$. In the case of  the complex simple Lie algebra of type $E_6$, denoted throughout this work by $\e6$, besides the split and the compact real forms, there are three more real forms, with signatures $-26$, $-14$ and $2$.

%%%%%%%%%%%%%%%%%%%%%%%%%%%%%%%%%%%%%%%%%%%%%%%%%%%%%%%%%%%%%%%%%%%%%%%%%%%%%%%%%%%%%%%%%%%%%%%%%%%%%%%%%%%%%%%%%%%%%%%%%%%%%%%%%%%%%%%
\section{Preliminaries about  Satake diagrams}\label{sec_preliminaresdeSatakes}

Let $\mathfrak{g}$ be an arbitrary semisimple Lie algebra over $\mathbb{R}$ and $k\colon \mathfrak{g}\times\mathfrak{g}\to\mathbb{R}$ its Killing form, which is nondegenerate. Recall \cite[Chapter III, \S7]{Helgasson} that a decomposition $\mathfrak{g}=\mathfrak{t}\oplus\mathfrak{p}$, for a subalgebra $\mathfrak{t}$ and a vector space  $\mathfrak{p}$, is called a \emph{Cartan decomposition} if there is a compact real form $\mathfrak{u}$ of $\mathfrak{g^\mathbb{C}}$ such that $\mathfrak{t}=\mathfrak{g}\,\cap\,\mathfrak{u}$ and $\mathfrak{p}=\mathfrak{g}\cap \imag\mathfrak{u}$. There always exists such a decomposition and any two Cartan decompositions are conjugate under an inner automorphism.
The automorphism $\theta\colon\mathfrak{g}\to\mathfrak{g}$ which sends $t+p$ to $t-p$ for any $t\in \mathfrak{t}$ and  $p\in \mathfrak{p}$,  is called a \emph{Cartan involution}.
In this case $k\vert_{\mathfrak{p}\times\mathfrak{p}}$ is positive definite and  $k\vert_{\mathfrak{t}\times\mathfrak{t}}$ is negative definite. Thus, it is said that $\mathfrak{t}$ is a \emph{maximal compactly imbedded subalgebra} of $\mathfrak{g}$. Observe that the signature of the Killing form coincides with $\dim\mathfrak{g}-2\dim\mathfrak{t}$.

Take any maximal abelian subspace $\mathfrak{a}$ of $\mathfrak{p}$. Its dimension is called  the \emph{real rank} of $\mathfrak{g}$ (which is independent of the choice of $\mathfrak{a}$).
For each $\lambda$ in the dual space $\mathfrak{ a}^*$ of $\mathfrak{ a}$, let
$\mathfrak{ g}_{\lambda}=\{x\in\mathfrak{ g}\mid [h,x]=\lambda(h)x\quad\forall h\in\mathfrak{ a}\}$. Then $\lambda$
is called a \emph{restricted root} if $\lambda\ne0$ and $\mathfrak{ g}_{\lambda}\ne0$. 
Denote by $\Sigma$ the set of restricted roots, which is an abstract root system (non necessarily reduced) according to \cite[Proposition~2.3.6]{parabolic},
and by $m_\lambda=\dim \mathfrak{g}_\lambda$ the \emph{multiplicity} of the restricted root.
Note that the simultaneous diagonalization of $\ad_\mathfrak{ g}\mathfrak{ a}$ gives the decomposition $\mathfrak{ g}=\mathfrak{g}_0\oplus\sum_{\lambda\in\Sigma}\mathfrak{g}_\lambda$, for $\mathfrak{g}_0=\mathfrak{a}\oplus \mathop{\rm Cent}_{\mathfrak{t}}(\mathfrak{a})$.
 
Now we combine the Cartan decomposition of a semisimple Lie algebra and the root space decomposition of its complexification.
Let $\mathfrak{h}$ be any maximal abelian subalgebra of $\mathfrak{g}$ containing $\mathfrak{a}$. Then $\mathfrak{h}$ is a Cartan subalgebra of $\mathfrak{g}$
(that is, $\mathfrak{h}^\mathbb{C}$ is a Cartan subalgebra of $\mathfrak{g}^\mathbb{C}$), and, if $\Delta$ denotes the root system of $\mathfrak{g}^\mathbb{C}$ relative to $\mathfrak{h}^\mathbb{C}$ and $\mathfrak{h}_\mathbb{R}:=\sum_{\alpha\in\Delta}\mathbb{R}h_\alpha$ (where  the element $h_\alpha\in[\mathfrak{g}^\mathbb{C}_\alpha,\mathfrak{g}^\mathbb{C}_{-\alpha}]$ is characterized by  $\alpha(h_\alpha)=2$), we get
$\mathfrak{h}_\mathbb{R}=\mathfrak{a}\oplus\imag(\mathfrak{h}\cap\mathfrak{t})$.
The restricted roots are exactly the nonzero restrictions of roots to $\mathfrak{a}\subset \mathfrak{h}^\mathbb{C} $, that is: If $\alpha\in \Delta$, denote by $\bar\alpha=\alpha\vert_{\mathfrak{a}}\colon\mathfrak{a}\to\mathbb{R}$. The roots  in $\Delta_0=\{\alpha\in\Delta\mid \bar \alpha=0\}$ are called the \emph{compact} roots and those in $\Delta\setminus\Delta_0$ the \emph{noncompact} roots.
Then $\Sigma=\{\bar \alpha\mid \alpha\in\Delta\setminus\Delta_0\}$. Note that $\alpha\in\Delta_0$ if and only if $\alpha(\mathfrak{h})\subset\imag\mathbb{R}$. Again $\Delta_0$ is an abstract root system on the Euclidean space of $\mathfrak{h}_\mathbb{R}$  spanned by its elements \cite[Proposition~2.3.8]{parabolic}.

Take a basis $B$  of the root system $\Delta$ adapted to our situation, that is, if we denote by $B_0=B\cap\Delta_0=\{\alpha\in B\mid \bar\alpha=0\}$, then the integral linear combinations of elements in $B_0$ with all the coefficients having the same sign coincide with the elements in $\Delta_0$ ($B_0$ is a basis of the root system $\Delta_0$). This   is equivalent to choose an ordering $\Delta^+$ such that for any $\alpha\in\Delta^+\setminus \Delta_0$ then $\sigma^*\alpha\in \Delta^+$, where $\sigma^*\alpha(h)=\overline{\alpha(\sigma(h))}$ for $\sigma$ the conjugation related to $\mathfrak{g}$.

The Satake diagram of the real algebra $\mathfrak{g}$ is defined as follows.
In the Dynkin diagram associated to the basis $B$, the roots in $B_0$ are denoted by a black circle $\bullet$ and the roots in $B\setminus B_0$ are denoted by a white circle $\circ$.
If $\alpha,\beta\in B\setminus B_0$ are such that $\bar \alpha=\bar \beta$, then $\alpha$ and $\beta$ are joined by a curved arrow.
Observe that the rank of $\mathfrak{g}^\mathbb{C}$ coincides with the real rank of $\mathfrak{g}$ plus the number of arrows on the Satake diagram plus the number of black nodes.
We enclose a list for each simple Lie algebra $\mathfrak{g}$ which is a real form of $\e6$.
The table contains the Dynkin diagram of $\bar B=\{\bar \alpha\mid \alpha\in B\setminus B_0\}$, which is a basis of the   root system $\Sigma$,
joint  with the multiplicities $m_\lambda$ and $m_{2\lambda}$  for $\lambda\in\bar B$. This table is extracted from \cite[Table VI]{Helgasson}, who obtained it from \cite{araki}.

A consistent extension of this notation is to define the Satake diagram of a compact semisimple Lie algebra to be the Dynkin diagram of the complexification $\mathfrak{g}^\mathbb{C}$   with all the nodes black.

\medskip

   \begin{center}{
 \begin{tabular}{ccccc}
  \hline\vspace{0.1cm}
  \vrule width 0pt height 14pt
  $ \mathfrak{g}$& Satake diagram of $(B,\theta)$ & Dynkin diagram of $\bar B$& $m_\lambda$ &$m_{2\lambda}$\cr
\hline\vspace{0.3cm}
EI&\begin{picture}(23,7)(4,-0.5)
\put(5,0){\circle{1}} \put(9,0){\circle{1}} \put(13,0){\circle{1}}
\put(17,0){\circle{1}} \put(21,0){\circle{1}}
\put(13,4){\circle{1}}
\put(5.5,0){\line(1,0){3}}
\put(9.5,0){\line(1,0){3}} \put(13.5,0){\line(1,0){3}} %\put(13,4.5){\line(0,1){3}}
\put(17.5,0){\line(1,0){3}}
\put(13,0.5){\line(0,1){3}}
\end{picture}&\begin{picture}(20,5)(4,-0.5)
\put(5,0){\circle{1}} \put(9,0){\circle{1}} \put(13,0){\circle{1}}
\put(17,0){\circle{1}} \put(21,0){\circle{1}}
\put(13,4){\circle{1}}
\put(5.5,0){\line(1,0){3}}
\put(9.5,0){\line(1,0){3}} \put(13.5,0){\line(1,0){3}} %\put(13,4.5){\line(0,1){3}}
\put(17.5,0){\line(1,0){3}}
\put(13,0.5){\line(0,1){3}}
\end{picture}&1&0\cr\vspace{0.3cm}
EII&\begin{picture}(23,5)(4,-0.5)
\put(5,0){\circle{1}} \put(9,0){\circle{1}} \put(13,0){\circle{1}}
\put(17,0){\circle{1}} \put(21,0){\circle{1}}
\put(13,4){\circle{1}}
\put(5.5,0){\line(1,0){3}}
\put(9.5,0){\line(1,0){3}} \put(13.5,0){\line(1,0){3}} %\put(13,4.5){\line(0,1){3}}
\put(17.5,0){\line(1,0){3}}
\put(13,0.5){\line(0,1){3}}
\put(4.5,1){$\scriptstyle \alpha_1$} \put(8.5,1){$\scriptstyle \alpha_3$}
\put(12.85,1){$\scriptstyle \alpha_4$} \put(16.5,1){$\scriptstyle \alpha_5$}
\put(20.5,1){$\scriptstyle \alpha_6$} \put(13.9,3.6){$\scriptstyle \alpha_2$}
\cbezier (6,-1)(11,-2.5)(15,-2.5)(20,-1)% cubic Bezier curve
 \put(20,-1){\vector(2,1){1}}
 \put(6,-1){\vector(-2,1){1}}
 \cbezier (10.2,-1)(12.5,-1.5)(14.5,-1.5)(16.0,-1)% cubic Bezier curve
 \put(16.0,-1){\vector(2,1){1}}
 \put(10.2,-1){\vector(-2,1){1}}
 \end{picture}&
\begin{picture}(25,5)(4,-0.5)
  \put(9,0){\circle{1}} \put(13,0){\circle{1}}
\put(17,0){\circle{1}} \put(21,0){\circle{1}}
\put(9.5,0){\line(1,0){3}} \put(13.5,0.3){\line(1,0){3}} \put(13.5,-0.1){\line(1,0){3}}
 \put(17.5,0){\line(1,0){3}} \put(14.3,-0.5){$<$}
  \put(8.5,-2){$\scriptstyle \alpha_1$}
\put(12.5,-2){$\scriptstyle \alpha_3$} \put(16.5,-2){$\scriptstyle \alpha_2$}
\put(20.5,-2){$\scriptstyle \alpha_4$}
\end{picture}
&$\left\{\begin{array}{ll}1 &(i=2,4)\\2&(i=1,3)\end{array}\right.$&$\begin{array}{l}0\\0\end{array}$\cr\vspace{0.3cm}
EIII&\begin{picture}(23,5)(4,-0.5)
\put(5,0){\circle{1}} \put(9,0){\circle*{1}} \put(13,0){\circle*{1}}
\put(17,0){\circle*{1}} \put(21,0){\circle{1}}
\put(13,4){\circle{1}}
\put(5.5,0){\line(1,0){3}}
\put(9.5,0){\line(1,0){3}} \put(13.5,0){\line(1,0){3}} %\put(13,4.5){\line(0,1){3}}
\put(17.5,0){\line(1,0){3}}
\put(13,0.5){\line(0,1){3}}
\put(4.5,1){$\scriptstyle \alpha_1$} \put(8.5,1){$\scriptstyle \alpha_3$}
\put(12.85,1){$\scriptstyle \alpha_4$} \put(16.5,1){$\scriptstyle \alpha_5$}
\put(20.5,1){$\scriptstyle \alpha_6$} \put(13.9,3.6){$\scriptstyle \alpha_2$}
\cbezier (6,-1)(11,-2.5)(15,-2.5)(20,-1)% cubic Bezier curve
 \put(20,-1){\vector(2,1){1}}
 \put(6,-1){\vector(-2,1){1}}
\end{picture}&
\begin{picture}(25,5)(4,-0.5)
   \put(13,0){\circle{1}}
\put(17,0){\circle{1}}
  \put(13.5,0.3){\line(1,0){3}} \put(13.5,-0.1){\line(1,0){3}}
  \put(14.3,-0.5){$<$}
 \put(12.5,-2){$\scriptstyle \alpha_1$} \put(16.5,-2){$\scriptstyle \alpha_2$}
\end{picture}
&$\left\{\begin{array}{ll}8 &(i=1)\\6&(i=2)\end{array}\right.$&$\begin{array}{l}1\\0\end{array}$\cr
EIV&\begin{picture}(23,5)(4,-0.5)
\put(5,0){\circle{1}} \put(9,0){\circle*{1}} \put(13,0){\circle*{1}}
\put(17,0){\circle*{1}} \put(21,0){\circle{1}}
\put(13,4){\circle*{1}}
\put(5.5,0){\line(1,0){3}}
\put(9.5,0){\line(1,0){3}} \put(13.5,0){\line(1,0){3}} %\put(13,4.5){\line(0,1){3}}
\put(17.5,0){\line(1,0){3}}
\put(13,0.5){\line(0,1){3}}
\put(4.5,1){$\scriptstyle \alpha_1$}
\put(20.5,1){$\scriptstyle \alpha_6$}
\end{picture}&\begin{picture}(25,5)(4,-0.5)
   \put(13,0){\circle{1}}
\put(17,0){\circle{1}}
  \put(13.5,0){\line(1,0){3}}
 \put(12.5,-2){$\scriptstyle \alpha_1$} \put(16.5,-2){$\scriptstyle \alpha_6$}
\end{picture}&8&0\cr&&&&\cr
\hline\vspace{0.001cm}
  \end{tabular}}

  Table I
  \end{center}

%%%%%%%%%%%%%%%%%%%%%%%%%%%%%%%%%%%%%%%%%%%%%%%%%%%%%%%%%%%%%%%%%%%%%%%%%%%%%%%%%%%%%%%%%%%%%%%%%%%%%%%%%%%%%%%%%%%%%%%%%%%%%%%%%%%%%%%
\section{Preliminaries about related algebras} %%Incluye simétricas

\subsection{Composition algebras}

A real \emph{composition algebra} $(C,n)$ is an $ \mathbb{R}$-algebra $C$ endowed with a nondegenerate quadratic form (the \emph{norm}) $n\colon  C \rightarrow \mathbb{R}$ which is multiplicative, that is,
$n(xy)=n(x)n(y)$ for all $x,y\in C$.
 Denote also by $n$ to the  {polar form}
$n(x,y):=n(x+y)-n(x)-n(y)$. A composition algebra is called \emph{split} if its norm is isotropic.

The unital composition algebras are called \emph{Hurwitz algebras}. Each Hurwitz algebra satisfies a quadratic equation
$$
x^2 - t(x)x + n(x)1=0,
$$
where the linear map $t(x):=n(x,1)$ is called the \emph{trace}. Besides it has a \emph{standard involution}   defined by
$\bar{x}:=t(x)1-x$, so that $n(x)=x\bar{x}$.

There is a standard process to construct these algebras starting from algebras of lower dimension, the so-called \emph{Cayley-Dickson doubling process}.
Let $A$ be a Hurwitz algebra with norm $n$ and let be $0\neq \alpha\in \mathbb{R}$.  Then the product $A\times A$ is endowed  with the following product:
\begin{equation}\label{eq_CayleyDicksonprocess}
(a,b)(c,d)=(ac+\alpha \bar{d}b, da+b\bar{c}).
\end{equation}
This new algebra, denoted by $CD(A,\alpha)$, has $(1,0)$ as a unit, it contains a copy of $A$ ($\equiv\{(a,0)\mid a\in A\}$)
and the element $u=(0,1)$ satisfies $u^2=\alpha1$, so that it can be identified with $A\oplus Au$. In particular $\dim (CD(A,\alpha))=2\dim A$.
Moreover   $CD(A,\alpha)$  is endowed with  the quadratic form
$
n((a,b))=n(a)-\alpha n(b)
$, being a  Hurwitz algebra if and only if $A$ is associative.

It can be easily proved that there are seven real Hurwitz algebras up to isomorphism, of dimensions  $1$, $2$, $4$ and $8$, namely:
\begin{itemize}
\item the ground field $\mathbb{R}$ (the involution is the identity);
\item the split algebra $\mathbb{R}\oplus \mathbb{R} \cong CD(\mathbb{R},1)$ (with   componentwise product and the exchange involution);
 \item the algebra of complex numbers $ \mathbb{C} = \span{1, i }\cong CD( \mathbb{R},-1)$ (the involution is the conjugation);
\item the matrix algebra $\Mat_{2\times2}(\mathbb{R})\cong CD(\mathbb{C},1)$ (the norm is given by the determinant);
 \item the quaternion division algebra   $ \mathbb{H}= \span{1, i, j, k} \cong CD( \mathbb{C},-1)$,   with $i^2=j^2=k^2=ijk=-1$ (the fundamental formula discovered by Hamilton in 1843);
\item    the   octonion division algebra $ \mathbb{O}= \span{1, i, j, k,l,il,jl,kl}\cong CD(\mathbb{H},-1)$, where the multiplication table is obtained from Equation~(\ref{eq_CayleyDicksonprocess}) for $CD(\mathbb{H},-1)=\mathbb{H}\oplus \mathbb{H}l$ with $l^2=-1$;
\item and
 the split octonion algebra $ \mathbb{O}_s\cong CD(\mathbb{H},1)$.
 This algebra has a \emph{standard basis}
 $\{e_1, e_2, u_1, u_2, u_3, v_1, v_2, v_3\} $ where $e_1$ and $e_2$ are orthogonal idempotents ($e_1+e_2=1$),
\begin{equation}\label{eq_baseestandar}
\begin{array}{lll}
e_1u_j=u_j=u_je_2,\ &  u_iv_i=-e_1, \    &    u_iu_{i+1}=v_{i+2}=-u_{i+1}u_i,       \\
e_2v_j=v_j=v_je_1,&   v_iu_i=-e_2,   &      v_iv_{i+1}=u_{i+2}=-v_{i+1}v_i,
\end{array}
\end{equation}
all the remaining products of basis elements are $0$,
and the   polar form of the norm   of two basis elements is zero except for $n(e_1,e_2)=1=n(u_i,v_i)$, $i=1,2,3$.
 \end{itemize}

\subsection{Symmetric composition algebras}

A  real composition algebra  $(C,n)$ is called a \emph{symmetric composition algebra}  if the (nondegenerate multiplicative)  quadratic form   satisfies
$n(x y,z)=n(x,y z)$
for any $x,y,z\in C$. The multiplication is usually denoted by $*$ instead of by juxtaposition, specially if there is some ambiguity.\smallskip

Again their dimensions are forced to be 1, 2, 4 or 8.  The only  examples are para-Hurwitz and pseudo-octonion algebras.

\begin{itemize}
\item
If $C$ is a Hurwitz algebra, the same vector space with new product
$$
x*y=\bar x\bar y
$$
for any $x,y\in C$ (and the same norm) is a symmetric composition algebra  called the \emph{para-Hurwitz algebra} attached to the Hurwitz algebra $C$. We will denote it by $pC$.

\item
Consider the algebra of $3\times3$ traceless matrices $\mathfrak{sl}(3,\mathbb{C})$ with the involution $x^*=p\bar x^t p^{-1}$ given by certain regular matrix $p$.
Note that  the (real) subspace of the antihermitian   matrices
$ S=\{x\in \mathfrak{sl}(3,\mathbb{C}) \mid x^* =  -x\}$
is closed for the product
\begin{equation*}\label{eq_productoOkubo}
x*y=\omega xy-\omega^2 yx-\frac{\omega-\omega^2}{3}\tr(xy)I_3,
\end{equation*}
 where   $\omega=e^{\frac{2\pi \imag}{3}}\in \mathbb{C}$ is a primitive cubic root of $1$ (so $\bar \omega=\omega^2$) and $I_3$ denotes the identity matrix.
 The obtained real algebra $(S,* )$, endowed with the norm
  $n(x) =-\frac12\tr(x^2)$, turns out to be a symmetric composition algebra, called \emph{pseudo-octonion} algebra or Okubo algebra.
  We will use the notation $\Ok$ in case
  $p=I_3$
and $\Ok_s$ when $p={\tiny \left( \begin{array}{ccc}
1 & 0 & 0\\
0 & 0 & 1\\
0 & 1 & 0
\end{array} \right)} $.
 Thus $\Ok$
coincides with $\mathfrak{su}(3)$ as a vector space and $\Ok_s$
   with $\mathfrak{su}(2,1)$. Observe that the norm $n\colon \Ok\to \mathbb{R}$ is definite, while  $n\colon \Ok_s\to \mathbb{R}$ is isotropic.
\end{itemize}

These pseudo-octonion algebras  are not
  isomorphic to any para-Hurwitz algebra. They were introduced by Susumu Okubo (see \cite[(4.35) and (4.9)]{libroOkubo}, where the definition appears with a minor modification) while he was working in Particle Physics. In particular, there are 4 real symmetric composition algebras of dimension 8: $p\mathbb{O}$ and  $p\mathbb{O}_s$ (also called para-octonion algebras), $ \Ok$ and $ \Ok_s$.

\subsection{Constructions of exceptional Lie algebras based on symmetric composition algebras}\label{subsec_constrcucciones}

Let $(S,\ast, q)$ be a symmetric composition algebra and let
$$
\mathfrak{o}(S,q)=\{d\in \End(S)\mid q(d(x),y)+q(x,d(y))=0\ \forall x,y\in S\}
$$
be the corresponding orthogonal Lie algebra. Consider the subalgebra of $\mathfrak{o}(S,q)^3$ defined by
$$
\mathfrak{tri}(S,\ast,q)=\{(d_0,d_1,d_2)\in \mathfrak{o}(S,q)^3 \mid
d_0(x\ast y)=d_1(x)\ast y + x\ast d_2(y)\ \forall x,y\in S\},
$$
 which is called the \emph{triality Lie algebra} in \cite{cuadradomagico}.
The order three automorphism $\vartheta$ given by
$$
\vartheta\colon  \mathfrak{tri}(S,\ast,q)\longrightarrow \mathfrak{tri}(S,\ast,q), \quad (d_0, d_1, d_2)\longmapsto (d_2, d_0, d_1),
$$
 is called the \emph{triality automorphism}. Take the element of $\mathfrak{tri}(S,\ast,q)$ (denoted by $\mathfrak{tri}(S)$ when there is no ambiguity) given by
\begin{equation}\label{eq_lostes}
t_{x,y}:=\left(\sigma_{x,y},\frac{1}{2}q(x,y)id-r_x l_y,\frac{1}{2}q(x,y)id-l_x r_y\right),
\end{equation}
where $\sigma_{x,y}(z)=q(x,z)y-q(y,z)x$, $r_x(z)=z\ast x$, and $l_x(z)=x\ast z$ for any $x,y,z\in S$.\smallskip

Let $(S,\ast,q)$ and $(S',\star ,q')$ be two real symmetric composition algebras.
Take nonzero scalars $\varepsilon_0,\varepsilon_1,\varepsilon_2\in\mathbb{R}$ and consider
$\varepsilon=(\varepsilon_0,\varepsilon_1,\varepsilon_2)$ and  the following vector space,
$$
\mathfrak{g}_\varepsilon(S,S') := \mathfrak{tri}(S,\ast, q)\oplus \mathfrak{tri}(S',\star, q')
\oplus (\bigoplus_{i=0}^2 \iota_i(S\otimes S')),
$$
where $\iota_i(S\otimes S')$ is just a copy of $S\otimes S'$ ($i=0,1,2$), and the anticommutative product
on $\mathfrak{g}_\varepsilon(S,S')$ is determined by the following conditions:
\begin{itemize}
\item   $\mathfrak{tri}(S,\ast, q)\oplus \mathfrak{tri}(S',\star, q')$ is a Lie subalgebra of $\mathfrak{ g}_\varepsilon(S,S')$,
\item $[(d_0,d_1,d_2), \iota_i(x\otimes x')]=\iota_i(d_i(x)\otimes x')$,

\item   $[(d'_0,d'_1,d'_2), \iota_i(x\otimes x')]=\iota_i(x\otimes d'_i(x'))$,
\item    $[\iota_i(x\otimes x'), \iota_{i+1}(y\otimes y')]=\varepsilon_{i+2}\iota_{i+2}((x*y)\otimes(x'\star y'))$,
\item   $[\iota_i(x\otimes x'),\iota_i(y\otimes y')] =\varepsilon_{i+1}\varepsilon_{i+2}( q'(x',y')\vartheta^i(t_{x,y}) + q(x,y)\vartheta'^i(t'_{x',y'}))$,
   \end{itemize}
for any $(d_0,d_1, d_2)\in \mathfrak{tri}(S)$, $(d'_0,d'_1, d'_2)\in \mathfrak{tri}(S')$,
   $x,y\in S$, $x',y'\in S'$,   $i=0,1,2$ indices taken modulo 3,
  and where  $\vartheta$ and $\vartheta'$  denote the corresponding triality automorphisms.

  When one of the involved symmetric composition algebras has dimension 2 (respectively 1) and the other one has dimension 8, the anticommutative algebra $\mathfrak{g}_\varepsilon(S,S')$ defined in this way turns out to be  a real  form of $\mathfrak{e}_6$ (respectively of  $\mathfrak{f}_4$)
  according to the following table (note that $\mathbb{R}=p\mathbb{R}$):

   \begin{center}{\begin{tabular}{lcr}
 \begin{tabular}{c|ccc|}
$  (1,1,1)$&$ \mathbb{R}$&$p\mathbb{C}$ & $p(\mathbb{R}+\mathbb{R})$
  \cr
 \hline
 $p\mathbb{O}, \,\Ok$&$\mathfrak{f}_{4, -52}$&$\mathfrak{e}_{6,-78}$& $\mathfrak{e}_{6,-26}$  \cr
 $p\mathbb{O}_s,\,\Ok_s$&$\mathfrak{f}_{4, 4}$&$\mathfrak{e}_{6,2}$& $\mathfrak{e}_{6,6}$  \cr
\hline
  \end{tabular} &\  &\begin{tabular}{c|ccc|}
$  (1,-1,1)$&$ \mathbb{R}$&$p\mathbb{C}$ & $p(\mathbb{R}+\mathbb{R})$
  \cr
 \hline
 $p\mathbb{O}, \,\Ok$&$\mathfrak{f}_{4, -20}$&$\mathfrak{e}_{6,-14}$& $\mathfrak{e}_{6,-26}$  \cr
 $p\mathbb{O}_s,\,\Ok_s$&$\mathfrak{f}_{4, 4}$&$\mathfrak{e}_{6,2}$& $\mathfrak{e}_{6,6}$  \cr
\hline
  \end{tabular}\cr\end{tabular}}
  \end{center}\vspace{0.5cm}

In other words, all the real forms of $\mathfrak{e}_6$ (and of  $\mathfrak{f}_4 $) can be obtained with this construction (details in \cite{algoreales}). Furthermore,  {all} the real forms of \textit{any} exceptional simple Lie algebra appear by choosing symmetric composition algebras $S$ and $S'$ of various dimensions.

\subsection{Jordan algebras}

A \emph{Jordan algebra} is a commutative (nonassociative) algebra satisfying the Jordan identity
$
(x^2y)x=x^2(yx).
$
An important example is
  $ H_3( \mathbb{O})=\{x=(x_{ij})\in \Mat_{3\times 3}( \mathbb{O})\mid x=\bar x^t\}$ with the product
  $x\cdot y=\frac12(xy+yx)$, which is a exceptional Jordan algebra   denoted  by $\mathbb{A}$. We will make use of the subspace of zero trace elements
  (for $\tr(x)=\sum_{i}x_{ii}\in\mathbb{R}$), which will be denoted by $\mathbb{A}_0$.

The real forms of the complexified algebra  $\mathbb{A}^\mathbb{C}$ are called real \emph{Albert algebras}, and they are quite involved in the real forms of the exceptional Lie algebras. There are three of them up to isomorphism: $\mathbb{A}$, $ H_3( \mathbb{O}_s)$, and a third one constructed in   \cite{Jacobsondeexcepcionales}
by doing slight  modifications to the construction of  $ H_3( \mathbb{O})$ (more precisely, considering  $x^*=p\bar x^tp^{-1}$ for $p={\tiny \left( \begin{array}{ccc}
1 & 0 & 0\\
0 & 0 & 1\\
0 & 1 & 0
\end{array} \right)}$ instead of being $p$ the identity matrix). We follow the approach  in \cite{reales} for a homogeneous description convenient for our aims.

Given a    real symmetric composition algebra $(S,*,
q)$ of dimension 8, and a set of
three nonzero scalars $\varepsilon=(\varepsilon_0,\varepsilon_1,\varepsilon_2)\in\{\pm1\}^3$, we define the commutative algebra
  $$
  {\mathbb A}_{( \varepsilon_0,\varepsilon_1,\varepsilon_2)}(S):=
\R^3 \oplus (\bigoplus\limits_{i=0}^2 \iota_i(S)),
$$
 where $\iota_i(S)$ is just a copy
of $S$ (i=0,1,2) and the  product is
given by
\begin{itemize}
 \item $(\alpha_0,\alpha_1,\alpha_2) (\beta_0,\beta_1,\beta_2)=(\alpha_0 \beta_0,\alpha_1 \beta_1,\alpha_2 \beta_2)$,
 \item
 $(\alpha_0,\alpha_1,\alpha_2)\iota_i(x)=\frac{1}{2}(\alpha_{i+1}+\alpha_{i+2})\iota_{i}(x),$
 \item $\iota_i(x)\iota_{i+1}(y)=\varepsilon_{i+2}\iota_{i+2}(x*y)$,
\item $\iota_i(x)\iota_{i}(y)=
    2\,\varepsilon_{i+1}\varepsilon_{i+2}q(x,y) (E_{i+1}+E_{i+2})$,
\end{itemize}
where the indices have been taken modulo 3 and $\{E_0,E_1,E_2\}$ denotes the
canonical basis of $\R^3$.
This algebra ${\mathbb A}_\varepsilon(S)$ is an Albert algebra, and conversely, the three real Albert algebras appear in this way (for suitable $\varepsilon$ and $S$).

This construction  of the   Albert algebras is related to the construction of the real forms of $ \mathfrak{f}_4$ in the above subsection. Namely,
there is a   Lie algebra isomorphism between the Lie algebra of derivations ${\rm Der}({\mathbb
A}_{\varepsilon}(S))$ and $ {\frak g}_{\varepsilon }(S, \R)= \mathfrak{tri}(S)\oplus (\oplus_{i=0}^2 \iota_i(S))$ (where $S\otimes \mathbb{R}$ has of course been identified with $S$)
given by the map
\begin{equation}\label{eq_themaprho}
\rho\colon {\frak g}_{\varepsilon }(S, \R)\to {\rm Der}({\mathbb
A}_{\varepsilon}(S)),
\end{equation}
where if $(d_0,d_1,d_2) \in {\frak{tri}}(S, *, q)$, then $\rho(d_0,d_1,d_2)$
is the derivation of ${\mathbb
A}_{\varepsilon}(S)$ given by $\rho(d_0,d_1,d_2)(\alpha_0, \alpha_1, \alpha_2)=0 $ and $
    \rho(d_0,d_1,d_2)(\iota_i(x))=\iota_i(d_i(x))$, and
    $\rho(\iota_i(x))$ is the derivation of ${\mathbb
A}_{\varepsilon}(S)$ given by  $2[l_{\iota_i(x)},
 l_{E_{i+1}}]$ for $l$ the left multiplication operator.

 %%%%%%%%%%%%%%%%%%%%%%%%%%%%%%%%%%%%%%%%%%%%%%%%%%%%%%%%%%%%%%%%%%%%%%%%%%%%%%%%%%%%%%%%%%%%%%%%%%%%%%%%%%%%%%%%%%%%%%%%%%%%%%%%%%%%%%%
\section{Satake diagram and facts on $\sig{-26}$}\label{sec26}

Consider  $ \mathfrak{g}= \der(\mathbb{A})\oplus\mathbb{A}_0$, which is endowed with a Lie algebra structure for the Lie bracket
$[x,y]:=[l_x,l_y]\in\der(\mathbb{A})$  if $x,y\in \mathbb{A}_0$ ($l$ the left multiplication operator)  and the natural action of $\der(\mathbb{A})$  on $\mathbb{A}_0$.
It is well known \cite{Jacobsondeexcepcionales} that $\mathfrak{g}$ is a real form of $ \mathfrak{e}_{6 }$ of signature $-26$, but also will be a consequence of the computation of its Satake diagram.

First, the decomposition $\mathfrak{g}=\mathfrak{t}\oplus \mathfrak{p}$, for
$$
\begin{array}{l}
  \mathfrak{t}= \der(\mathbb{A}),  \\
 \mathfrak{p}= \mathbb{A}_0,
\end{array}
$$
is a Cartan decomposition, since
  $k\vert_{ \mathfrak{p}}$ is positive definite and $k\vert_{ \mathfrak{t}}$ is negative definite, being $k\colon \mathfrak{g}\times  \mathfrak{g}\to\mathbb{R}$ the Killing form.  (In particular, the signature of $\mathfrak{g}$ is $-52 +26 =-26$.) The  {maximal compactly imbedded subalgebra}   $\mathfrak{t}$ is of type $F_4$.

For unifying the notation, note the isomorphism   $\mathbb{A}_{(1,1,1)}(p\mathbb{O})\cong \mathbb{A}$ by means of
\begin{equation*}\label{eq_isomorfismodealbert}
(\alpha_0,\alpha_1,\alpha_2)+\iota_0(x)+\iota_1(y)+\iota_2(z)\mapsto\left(\begin{array}{ccc}
\alpha_0&2z&2\bar y\\ 2\bar  z&\alpha_1&2x\\ 2y&2\bar x& \alpha_2\end{array}\right).
\end{equation*}
So, we will work with $\mathfrak{g}\equiv\rho(\mathfrak{g}_{(1,1,1)}(p\mathbb{O},\mathbb{R}) )\oplus (\mathbb{A}_{(1,1,1)}(p\mathbb{O}))_0$,
being $\rho$ the map defined in Equation~(\ref{eq_themaprho}).

Second, we find $ \mathfrak{h}=\mathfrak{ a}\oplus  (\mathfrak{h}\cap \mathfrak{t})$ a Cartan subalgebra of $ \mathfrak{g}$ such that $ \mathfrak{a}$ is a maximal abelian subalgebra of $ \mathfrak{p}$. If we fix the basis of the octonion division algebra
\begin{equation}\label{eq_basedeoctondivision}
\{e_0,e_1,e_2,e_3,e_4,e_5,e_6,e_7\}:=\{1,i,j,k,l,il,jl,kl\},
\end{equation}
(the bilinear form $\frac12q$ relative to this basis is the identity matrix of size 8), then we can take
$$
\begin{array}{l}
\ \mathfrak{a}=\langle\{E_0-E_1,E_2-E_0\}\rangle\subset\mathbb{A}_0,\\
\mathfrak{h}\cap \mathfrak{t} =\langle\{\rho(t_{e_0,e_1}),\rho(t_{e_2,e_3}),\rho(t_{e_4,e_5}),\rho(t_{e_6,e_7 })\}\rangle\subset \rho(\mathfrak{tri}(p\mathbb{O})),
\end{array}
$$
with the notation for the elements in the triality Lie algebra considered in Equation~(\ref{eq_lostes}).

Third, we decompose the complex Lie algebra $ \mathfrak{g}^\mathbb{C} $ relative to its Cartan algebra
$ \mathfrak{h}^\mathbb{C}= \mathfrak{h}\otimes_\mathbb{R}\mathbb{C}$. Take as a basis
$$
\begin{array}{lll}\vspace{0.1cm}
h_1=\frac12\rho(t_{e_0,e_1}), & h_3= \frac12\rho(t_{e_4,e_5}), & h_5= E_2-E_0, \\
h_2=\frac12\rho(t_{e_2,e_3}), &  h_4=\frac12 \rho(t_{e_6,e_7}), &h_6= E_0-E_1 ,
\end{array}
$$
and an arbitrary element $h=\sum_{i=1}^6w_ih_i\in  \mathfrak{g}^\mathbb{C}$.
A long but straightforward computation shows that $\ad h$ diagonalizes $ \mathfrak{g}^\mathbb{C}$ with set of eigenvalues $\Delta$ union of:
\begin{itemize}
\item $\{\imag(\pm w_j\pm w_k)\mid j\ne k, 1\le j,k\le 4\}$ in $\rho(\mathfrak{tri}(p\mathbb{O}))$,
\item $\{\pm\imag w_j\pm\frac{w_5+w_6}{2}\mid 1\le j\le 4\}$ in $\rho(\iota_0(p\mathbb{O}))\oplus \iota_0(p\mathbb{O})$,
\item $\{\frac12(\imag (\varepsilon_1 w_1+\varepsilon_2 w_2+\varepsilon_3 w_3+\varepsilon_4 w_4)\pm{(2w_5-w_6)})\mid \varepsilon_i=\pm1,\Pi_{i=1}^4\varepsilon_i=1\}$ in $\rho(\iota_1(p\mathbb{O}))\oplus \iota_1(p\mathbb{O})$,
\item $\{\frac12(\imag (\varepsilon_1 w_1+\varepsilon_2 w_2+\varepsilon_3 w_3+\varepsilon_4 w_4)\pm{(w_5-2w_6)})\mid \varepsilon_i=\pm1,\Pi_{i=1}^4\varepsilon_i=-1\}$ in $\rho(\iota_2(p\mathbb{O}))\oplus \iota_2(p\mathbb{O})$.
    \end{itemize}

Take $B=\{\alpha_1,\alpha_2,\alpha_3,\alpha_4,\alpha_5,\alpha_6\}\subset \Delta$ for $\alpha_i\colon  \mathfrak{g}^\mathbb{C}\to\mathbb{C}$ given by
$$
\begin{array}{l}
\alpha_1(h)= \frac12( \ii w_1-\ii w_2-\ii w_3-\ii w_4+w_5-2w_6)  ,\\
\alpha_2(h)= \ii (w_1+w_2),\\
\alpha_3(h)=\ii(-w_1+w_2),\\
\alpha_4(h)=\ii(-w_2+w_3),\\
\alpha_5(h)=\ii(-w_3+w_4),\\
\alpha_6(h)= \frac12(-2\ii w_4+w_5+w_6 ) .
\end{array}
$$
Again it is a long and direct computation to check that $B$ is a basis of the root system $\Delta$, since $\Delta\subset \sum_{i=1}^6\mathbb{Z}_{\ge0}\alpha_i\cup\sum_{i=1}^6\mathbb{Z}_{\le0}\alpha_i$.
%%%Surge la duda de si merece poner las raíces que todo el mudo se sabe o me lo ahorro
%
The choice of $B$ is well adapted to the situation, since
  the set of compact roots in the basis, $B_0=\{\alpha\in B\mid \alpha( \mathfrak{a})=0\}=\{\alpha\in B\mid \alpha(h_5)=\alpha(h_6)=0\}$, that is,
$$
B_0=\{\alpha_2,\alpha_3,\alpha_4,\alpha_5\},
$$
turns out to be a basis of the set of compact roots $\Delta_0=\{\imag(\pm w_j\pm w_k)\mid j\ne k, 1\le j,k\le 4\}$ (which is of type $D_4$).
As $\alpha_1(w_5h_5+w_6h_6)=\frac12( w_5-2w_6)\ne \frac12( w_5+w_6)=\alpha_6(w_5h_5+w_6h_6)$, then the nodes related to $\alpha_1$ and $\alpha_6$ are not joined in the Satake diagram, which is
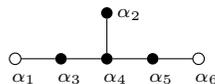
\begin{center}{
\begin{picture}(25,5)(4,-0.5)
\put(5,0){\circle{1}} \put(9,0){\circle*{1}} \put(13,0){\circle*{1}}
\put(17,0){\circle*{1}} \put(21,0){\circle{1}}
\put(13,4){\circle*{1}}
\put(5.5,0){\line(1,0){3}}
\put(9.5,0){\line(1,0){3}} \put(13.5,0){\line(1,0){3}} %\put(13,4.5){\line(0,1){3}}
\put(17.5,0){\line(1,0){3}}
\put(13,0.5){\line(0,1){3}}
\put(4.7,-2){$\scriptstyle \alpha_1$} \put(8.7,-2){$\scriptstyle \alpha_3$}
\put(12.7,-2){$\scriptstyle \alpha_4$} \put(16.7,-2){$\scriptstyle \alpha_5$}
\put(20.7,-2){$\scriptstyle \alpha_6$} \put(13.9,3.6){$\scriptstyle \alpha_2$}
\end{picture}
}\end{center}\vskip0.4cm

Finally the set of restricted roots is $\Sigma=\pm \{ \frac{   w_5 + w_6 }2,\frac{   w_5-2  w_6 }2,\frac{   2w_5 - w_6 }2\}$, which  is a root system of type $A_2$ (with basis $\bar B$), more precisely,
$\Sigma=\pm\{\overline{\alpha}_1,\overline{\alpha}_6,\overline{\alpha_1+\alpha_6} \}$.
The multiplicities are
$$m_{\overline{\alpha}_1}=8, \quad    m_{\overline{\alpha}_6 }=8.$$

%%%%%%%%%%%%%%%%%%%%%%%%%%%%%%%%%%%%%%%%%%%%%%%%%%%%%%%%%%%%%%%%%%%%%%%%%%%%%%%%%%%%%%%%%%%%%%%%%%%%%%%%%%%%%%%%%%%%%%%%%%%%%%%%%%%%%%%
\section{Satake diagram and facts on $\sig{-14}$}\label{sec14}

Consider  $ \mathfrak{g}= \mathfrak{g}_{(1,-1,1)}(p\mathbb{O},p\mathbb{C})$, which is a real form of $ \mathfrak{e}_{6 }$ of signature $-14$ according to   Section~\ref{subsec_constrcucciones}, but also will be a consequence of the next computations.

Take the decomposition $\mathfrak{g}=\mathfrak{t}\oplus \mathfrak{p}$, for
$$
\begin{array}{l}
  \mathfrak{t}= \mathfrak{tri}(p\mathbb{O})\oplus \mathfrak{tri}(p\mathbb{C})\oplus \iota_1( p\mathbb{O}\otimes p\mathbb{C}), \\
 \mathfrak{p}=  \iota_0( p\mathbb{O}\otimes p\mathbb{C})\oplus \iota_2( p\mathbb{O}\otimes p\mathbb{C}),
\end{array}
$$
which is a Cartan decomposition  since $\mathfrak{t}\oplus\imag\mathfrak{p}\cong \mathfrak{g}_{(1, 1,1)}(p\mathbb{O},p\mathbb{C}) \cong\sig{-78}$ is compact. Thus
  $k\vert_{ \mathfrak{p}}$ is positive definite and $k\vert_{ \mathfrak{t}}$ is negative definite, being $k\colon \mathfrak{g}\times  \mathfrak{g}\to\mathbb{R}$ the Killing form.
  (In particular, the signature of $\mathfrak{g}$ is $-46 +32 =-14$.)
  The  {maximal compactly imbedded subalgebra}   $\mathfrak{t}$ is of type $D_5+Z$.

Second, we find $ \mathfrak{h}=\mathfrak{ a}\oplus  (\mathfrak{h}\cap \mathfrak{t})$ a Cartan subalgebra of $ \mathfrak{g}$ such that $ \mathfrak{a}$ is a maximal abelian subalgebra of $ \mathfrak{p}$. With the basis of  $p\mathbb{O}$  fixed in Equation~(\ref{eq_basedeoctondivision}), in which the two first elements $\{e_0,e_1\}$ can be considered as a basis of the paracomplex algebra $p\mathbb{C}$, we take
$$
\begin{array}{l}
\ \mathfrak{a}=\langle\{\iota_0(e_0\otimes e_1),\iota_0(e_1\otimes e_0)\}\rangle,\\
\mathfrak{h}\cap \mathfrak{t} =\langle\{t_{e_2,e_3},t_{e_4,e_5},t_{e_6,e_7 },(0,\sigma'_{e_0,e_1},-\sigma'_{e_0,e_1})\}\rangle,
\end{array}
$$
with the notations $t_{x,y}$ given in Equation~(\ref{eq_lostes}), and where the primes are again used for the second symmetric composition involved, in this case, $p\mathbb{C}$.

Now we decompose the complex Lie algebra $ \mathfrak{g}^\mathbb{C} $ relative to its Cartan algebra
$ \mathfrak{h}^\mathbb{C}= \mathfrak{h}\otimes_\mathbb{R}\mathbb{C}$. Take as a basis
$$
\begin{array}{lll}\vspace{0.1cm}
h_1=\frac12 t_{e_2,e_3}, & h_3=\frac12  t_{e_6,e_7}, & h_5=\frac12  \iota_0(e_0\otimes e_1),\\
h_2=\frac12 t_{e_4,e_5}, &  h_4=\frac14 (0,\sigma'_{e_0,e_1},-\sigma'_{e_0,e_1}), &h_6=\frac12  \iota_0(e_1\otimes e_0),
\end{array}
$$
and an arbitrary element $h=\sum_{i=1}^6w_ih_i\in  \mathfrak{g}^\mathbb{C}$.
We get that  $\ad h$ acts on
\begin{itemize}
\item the eigenvector $t_{e_0,e_2}+\imag t_{e_0,e_3}+\iota_0((e_2+\imag e_3)\otimes e_1) $ with eigenvalue $-\imag w_1-w_5 $,
\item the eigenvector $t_{e_2,e_4}+\imag t_{e_3,e_4}+\imag t_{e_2,e_5}-t_{e_3,e_5} $ with eigenvalue $-\imag w_1-\imag w_2 $,
\item the eigenvector $\iota_0(e_0\otimes e_0)-\iota_0(e_1\otimes e_1)+t_{e_0,e_1}+t'_{e_0,e_1} $ with eigenvalue $w_5+w_6 $,
\end{itemize}
where that kind of eigenvectors spans $\mathfrak{tri}(p\mathbb{O})\oplus \mathfrak{tri}(p\mathbb{C})\oplus \iota_0( p\mathbb{O}\otimes p\mathbb{C})$, and on
\begin{itemize}
\item the eigenvector $\iota_1((e_0+\imag e_1)\otimes (e_0+\imag e_1) ) +\iota_2((e_0-\imag e_1)\otimes (\imag e_0+e_1) )$ with eigenvalue $\frac12( \imag w_1+\imag w_2+\imag w_3-\imag w_4-w_5-w_6) $.
\end{itemize}
Analogously, it is checked that $\ad h$ diagonalizes $ \mathfrak{g}^\mathbb{C}$ with eigenvalues
$$
\begin{array}{rl}
 \Delta=&\{\pm\ii w_k\pm w_5,\pm\ii w_k\pm w_6\mid k=1,2,3\}\\
&\cup\,\{\ii(\pm w_1\pm w_2 ),\ii(\pm w_1\pm w_3 ),\ii(\pm w_2\pm w_3 ) \}\\
&\cup\,\{\pm w_5\pm w_6\}\\
&\cup\,\{\frac12(\ii(\pm w_1\pm w_2\pm w_3\pm w_4)\pm w_5\pm w_6)\}.
\end{array}
$$
The choice of a basis is usually the most difficult task. In this case, a suitably adapted  election is
  $B=\{\alpha_1,\alpha_2,\alpha_3,\alpha_4,\alpha_5,\alpha_6\}\subset \Delta$ for $\alpha_i\colon  \mathfrak{g}^\mathbb{C}\to\mathbb{C}$ given by
$$
\begin{array}{l}
\alpha_1(h)= \frac12(-\ii w_1-\ii w_2-\ii w_3-\ii w_4-w_5+w_6)  ,\\
\alpha_2(h)=-\ii w_1+w_5,\\
\alpha_3(h)=\ii(w_2+w_3),\\
\alpha_4(h)=\ii(w_1-w_2),\\
\alpha_5(h)=\ii(w_2-w_3),\\
\alpha_6(h)= \frac12(-\ii w_1-\ii w_2+\ii w_3+\ii w_4-w_5+w_6 ) ,
\end{array}
$$
because it is straightforward to check   $\Delta\subset \sum_{i=1}^6\mathbb{Z}_{\ge0}\alpha_i\cup\sum_{i=1}^6\mathbb{Z}_{\le0}\alpha_i$
(so that $B$ is a basis of the root system $\Delta$)
and also   $B_0=\{\alpha\in B\mid \alpha( \mathfrak{a})=0\}=\{\alpha\in B\mid \alpha(h_5)=\alpha(h_6)=0\}$ is
$$
B_0=\{\alpha_3,\alpha_4,\alpha_5\} ,
$$
 a basis of $\Delta_0=\{\ii(\pm w_1\pm w_2 ),\ii(\pm w_1\pm w_3 ),\ii(\pm w_2\pm w_3 ) \}$ (a root system of type $A_3$).
As $\alpha_1(w_5h_5+w_6h_6)=\frac12(-w_5+w_6)=\alpha_6(w_5h_5+w_6h_6)$,   the nodes related to $\alpha_1$ and $\alpha_6$ are joined by an arrow and then the Satake diagram is
\vskip0.3cm

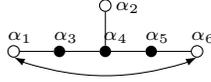
\begin{center}{
\begin{picture}(23,5)(4,-0.5)
\put(5,0){\circle{1}} \put(9,0){\circle*{1}} \put(13,0){\circle*{1}}
\put(17,0){\circle*{1}} \put(21,0){\circle{1}}
\put(13,4){\circle{1}}
\put(5.5,0){\line(1,0){3}}
\put(9.5,0){\line(1,0){3}} \put(13.5,0){\line(1,0){3}} %\put(13,4.5){\line(0,1){3}}
\put(17.5,0){\line(1,0){3}}
\put(13,0.5){\line(0,1){3}}
\put(4.5,1){$\scriptstyle \alpha_1$} \put(8.5,1){$\scriptstyle \alpha_3$}
\put(12.85,1){$\scriptstyle \alpha_4$} \put(16.5,1){$\scriptstyle \alpha_5$}
\put(20.5,1){$\scriptstyle \alpha_6$} \put(13.9,3.6){$\scriptstyle \alpha_2$}
\cbezier (6,-1)(11,-2.5)(15,-2.5)(20,-1)% cubic Bezier curve
 \put(20,-1){\vector(2,1){1}}
 \put(6,-1){\vector(-2,1){1}}
\end{picture}
}\end{center}\vskip0.4cm

In this case the set of restricted roots is $\Sigma=\{\pm w_5,\pm w_6,\pm w_5\pm w_6,\frac12(\pm w_5\pm w_6)\}$,
which coincides with
$$
\Sigma=\pm\{\overline{\beta}_1,\overline{\beta}_2,\overline{\beta}_1+ \overline{\beta}_2,2\overline{\beta}_1+\overline{\beta}_2,2\overline{\beta}_1,2\overline{\beta}_1+2\overline{\beta}_2\}
$$
for $\beta_1=\alpha_1+\alpha_2+\alpha_3+\alpha_4$ and $\beta_2=-\alpha_2$, since
$\overline{\beta_1}=\frac{w_5+w_6}{2}$ and $\overline{\beta_2}=-w_5$. Thus $\Sigma$ is a nonreduced root system of type $BC_2$.
The multiplicities of the restricted roots are
$$
\begin{array}{l}
m_{\frac{w_5+w_6}{2}}=8,\\
m_{w_5+w_6}=1,\\
m_{-w_5}=6.
\end{array}
$$\smallskip

We can observe that the noncompact (obviously long) roots in the basis ($\alpha_1$ and $\alpha_2$) have (restricted) multiplicity different from 1 (so that their root vectors cannot be used to obtain $\mathbb{Z}$-gradings as in \cite[Theorem 3]{Cheng}, as recalled in Introduction), but the maximal root $\alpha_1+2\alpha_2+2\alpha_3+3\alpha_4+2\alpha_5+\alpha_6=w_5+w_6$ has (restricted) multiplicity 1 and the $\mathbb{Z}$-grading does appear.

%%%%%%%%%%%%%%%%%%%%%%%%%%%%%%%%%%%%%%%%%%%%%%%%%%%%%%%%%%%%%%%%%%%%%%%%%%%%%%%%%%%%%%%%%%%%%%%%%%%%%%%%%%%%%%%%%%%%%%%%%%%%%%%%%%%%%%%
\section{Satake diagram and facts on $\sig2$}\label{sec2}

Consider  $ \mathfrak{g}= \mathfrak{g}_{(1, 1,1)}(p\mathbb{O}_s,p\mathbb{C})$, which is a real form of $ \mathfrak{e}_{6 }$ of signature $2$ according to
 Section~\ref{subsec_constrcucciones}, although this fact will also be a consequence of the next computations. % of its Satake diagram. %(Also for $\alpha=(1,-1,1)$.)
The description of a Cartan decomposition is more involved than in the previous cases, so we will work with this algebra a little bit before giving one.
Recall that $d=\sigma'_{e_0,e_1}\in \mathfrak{o}(p\mathbb{C},n)$ acts in the paracomplex algebra by sending $e_0$ to $2e_1$ and $e_1$ to $-2e_0$, and
$\mathfrak{tri}(p\mathbb{C})=\{(\beta_0d,\beta_1d,\beta_2d)\mid \sum_{i=0}^2 \beta_i=0\}$ is a two-dimensional abelian algebra.
Take
$$\begin{array}{lll}\vspace{0.15cm}
h_1=t_{e_1 ,e_2 },\quad & h_3=t_{u_2 ,v_2 },\quad& h_5=\frac14(\sigma'_{e_0,e_1},\sigma'_{e_0,e_1},-2\sigma'_{e_0,e_1}),\\
 h_2=t_{u_1 ,v_1 },  &h_4=t_{u_3,v_3 },& h_6=\frac14(\sigma'_{e_0,e_1},-2\sigma'_{e_0,e_1},\sigma'_{e_0,e_1}),
\end{array}
$$
where we are taking the standard basis of $\mathbb{O}_s$  as in Equation~(\ref{eq_baseestandar}).
It turns out that $\sum_{i=1}^6\mathbb{C}h_i$ is a Cartan subalgebra of  $ \mathfrak{g}^\mathbb{C}$. More concretely, an arbitrary element
$ \sum_{i=1}^4 w_ih_i\in \mathfrak{h}_0=\langle h_1,h_2,h_3,h_4\rangle$    is ad-diagonalizable with eigenvalues
\begin{itemize}
\item $0$ in $ \mathfrak{tri}(p\mathbb{C})\oplus \mathfrak{h}_0$,
\item $\Phi_l=\{ \pm w_j\pm w_k \mid j\ne k, 1\le j,k\le 4\}$ in $ \mathfrak{tri}(p\mathbb{O}_s)$ (to be precise, $\Phi_l\cup\{0\}$),
\item $\Phi_0=\{\pm w_j\mid 1\le j\le 4\}$ in $\iota_0(p\mathbb{O}_s\otimes p\mathbb{C})$,
\item $\Phi_1=\{\frac12(\varepsilon_1 w_1+\varepsilon_2 w_2+\varepsilon_3 w_3+\varepsilon_4 w_4) \mid \varepsilon_i=\pm1,\Pi_{i=1}^4\varepsilon_i=1\}$
in $\iota_1(p\mathbb{O}_s\otimes p\mathbb{C}))$,
\item $\Phi_2=\{\frac12(\varepsilon_1 w_1+\varepsilon_2 w_2+\varepsilon_3 w_3+\varepsilon_4 w_4) \mid \varepsilon_i=\pm1,\Pi_{i=1}^4\varepsilon_i=-1\}$
in $\iota_2(p\mathbb{O}_s\otimes p\mathbb{C}))$;
\end{itemize}
and $w_5h_5+w_6h_6$ acts with eigenvalues
\begin{itemize}
\item $0$ in $ \mathfrak{tri}(p\mathbb{C})\oplus \mathfrak{tri}(p\mathbb{O}_s)$,
\item  $\mp\frac\imag2(w_5+w_6)$ in $\iota_0(x\otimes(e_0\pm\imag e_1))$ for all $x\in p\mathbb{O}_s$,
\item  $\pm\imag(\frac12 w_5-w_6)$ in $\iota_1(x\otimes(e_0\pm\imag e_1))\in\iota_1(p\mathbb{O}_s\otimes p\mathbb{C}))$,
\item  $\pm\imag(w_5-\frac12 w_6)$ in $\iota_2(x\otimes(e_0\pm\imag e_1))\in\iota_2(p\mathbb{O}_s\otimes p\mathbb{C}))$.
\end{itemize}
In particular $\mathfrak{g}_0=
\mathfrak{g}_{(1, 1,1)}(p\mathbb{O}_s, \mathbb{R}e_0)$, which  is a real split subalgebra of $\mathfrak{g}$ isomorphic to $\mathfrak{f}_{4,4}$,
has $\mathfrak{h}_0 $ as a Cartan subalgebra   with roots $\Phi=\Phi_l\cup\Phi_0\cup\Phi_1\cup\Phi_2$
(observe that $\Phi_l$ are the long roots of $\Phi$ and $\Phi_0\cup\Phi_1\cup\Phi_2$ are the short ones, which correspond to   $\sum_{i=0}^2\iota_i(p\mathbb{O}_s\otimes \mathbb{R}e_0)$).

As
$\{\frac12(w_1-w_2-w_3-w_4),w_4,w_3-w_4,w_2-w_3\}$ is a basis of the root system $\Phi$,
then we have an ordering with positive roots
$\Phi^+=\{w_j\pm w_k\mid j<k\}\cup\{w_j\mid 1\le j\le 4\}\cup\{\frac12( w_1+\varepsilon_2 w_2+\varepsilon_3 w_3+\varepsilon_4 w_4) \mid \varepsilon_i=\pm1\}$.
Take, for each $\alpha\in\Phi^+$, elements $e_\alpha\in (\mathfrak{g}_0)_\alpha$ and
$f_\alpha\in (\mathfrak{g}_0)_{-\alpha}$ such that $[e_\alpha,f_\alpha]=h_\alpha$ (defined  as in Section~\ref{sec_preliminaresdeSatakes} by $\alpha(h_\alpha)=2$).
That implies that $k_0(e_\alpha,f_\alpha)>0$ (and $k_0((\mathfrak{g}_0)_{ \alpha},(\mathfrak{g}_0)_{ \beta})=0$ if $\alpha+\beta\ne0$) for $k_0$ the Killing form of $\mathfrak{g}_0$. In particular
$\imag \mathfrak{h}_0\oplus\langle\{ e_\alpha-f_\alpha\mid \alpha\in \Phi^+\}\rangle\oplus  \langle \{\imag(e_\alpha+f_\alpha)\mid \alpha\in  \Phi^+\}\rangle\cong\mathfrak{f}_{4,-52}$
is a compact real form of $(\mathfrak{g}_0)^\mathbb{C}$.

Note that, for each $i=0,1,2$, the map $\Psi_i\colon \mathfrak{g}\to\mathfrak{g}$ given by
\begin{itemize}
\item $\Psi_i\vert_{\mathfrak{tri}(p\mathbb{O}_s)\oplus \mathfrak{tri}(p\mathbb{C})}=\id$,
\item $\Psi_i\vert_{\iota_i(p\mathbb{O}_s\otimes p\mathbb{C})}=-\id$,
\item $\Psi_i(\iota_j(x\otimes e_k))= (-1)^{k+1}\iota_j(x\otimes e_{k+1})$ if $j\ne i$, $x\in p\mathbb{O}_s$ and $k=0,1$ (mod~2),
\end{itemize}
is an automorphism. 
In particular $k(\Psi_i(e_\alpha),\Psi_i(f_\alpha))=k(e_\alpha,f_\alpha)$ (a positive multiple of $k_0(e_\alpha,f_\alpha)  $) for each $\alpha\in\Phi^+$.
Since $\iota_i(p\mathbb{O}_s\otimes  \mathbb{R}e_1)=\Psi_{i+1}(\iota_i(p\mathbb{O}_s\otimes  \mathbb{R}e_0))$, then
$$
\begin{array}{l}\vspace{0.07cm}
\mathfrak{tri}(p\mathbb{C})\oplus \imag \mathfrak{h}_0\oplus\langle\{ e_\alpha-f_\alpha\mid \alpha\in \Phi^+\}\rangle\oplus
\langle\{ \imag(e_\alpha+f_\alpha)\mid \alpha\in  \Phi^+\}\rangle\oplus\\
 \oplus_{i=0}^2 \left(\langle\{ \Psi_i(e_\alpha-f_\alpha)\mid \alpha\in \Phi^+\cap \Phi_{i+1}\}\rangle\oplus
\langle\{ \imag\Psi_i(e_\alpha+f_\alpha)\mid \alpha\in  \Phi^+\cap \Phi_{i+1}\}\rangle \right)
\end{array}
$$
is a compact real form of $\mathfrak{g}^\mathbb{C}$. That means that a Cartan decomposition of $\mathfrak{g}$ is the following:
$$
\begin{array}{l}\vspace{0.07cm}
\mathfrak{t}=\mathfrak{tri}(p\mathbb{C})\oplus \langle\{ e_\alpha-f_\alpha\mid \alpha\in \Phi^+\}\rangle\oplus\left(\oplus_{i=0}^2 \langle\{ \Psi_i(e_\alpha-f_\alpha)\mid \alpha\in \Phi^+\cap \Phi_{i+1}\}\rangle\right),\\
\mathfrak{p}=\mathfrak{h}_0\oplus  \langle\{  e_\alpha+f_\alpha\mid \alpha\in  \Phi^+\}\rangle\oplus \left(\oplus_{i=0}^2 \langle\{ \Psi_i(e_\alpha+f_\alpha)\mid \alpha\in \Phi^+\cap \Phi_{i+1}\}\rangle\right),
\end{array}
$$
  and the signature of $\mathfrak{g}$ turns out to be $38-36=2$.
  Thus $ \mathfrak{h}=\mathfrak{ a}\oplus  (\mathfrak{h}\cap \mathfrak{t})$ is a suitably adapted Cartan subalgebra for
$$
\begin{array}{l}
 \mathfrak{a}=\mathfrak{h}_0\subset \mathfrak{p},\\
\mathfrak{h}\cap \mathfrak{t} = \mathfrak{tri}(p\mathbb{C}),
\end{array}
$$
and we have already done the simultaneous diagonalization of the complex Lie algebra $ \mathfrak{g}^\mathbb{C} $ relative to
$ \mathfrak{h}^\mathbb{C}$:
An arbitrary element $ \sum_{i=1}^6w_ih_i\in  \mathfrak{h}^\mathbb{C}$
 diagonalizes $ \mathfrak{g}^\mathbb{C}$ with eigenvalues
 $$
 \begin{array}{ll}\vspace{0.1cm}
 \Delta=&
 \{ \pm w_j\pm w_k \mid j\ne k, 1\le j,k\le 4\}
 \cup\{\pm w_j\pm\frac{\imag}2(w_5+w_6)\mid 1\le j\le 4\}\\\vspace{0.1cm}
 &\cup\,\{\frac12(\varepsilon_1 w_1+\varepsilon_2 w_2+\varepsilon_3 w_3+\varepsilon_4 w_4)\pm\imag( \frac12 w_5-w_6) \mid \varepsilon_i=\pm1,\Pi_{i=1}^4\varepsilon_i=1\}\\
 &\cup\,\{\frac12(\varepsilon_1 w_1+\varepsilon_2 w_2+\varepsilon_3 w_3+\varepsilon_4 w_4) \pm\imag( w_5- \frac12 w_6)\mid \varepsilon_i=\pm1,\Pi_{i=1}^4\varepsilon_i=-1\}
 .
 \end{array}$$

 Take $B=\{\alpha_1,\alpha_2,\alpha_3,\alpha_4,\alpha_5,\alpha_6\}\subset \Delta$ for $\alpha_i\colon  \mathfrak{g}^\mathbb{C}\to\mathbb{C}$ given by
$$
\begin{array}{l}
\alpha_1(h)=\frac12(w_1-w_2-w_3-w_4+\imag(2w_5-w_6))    ,\\
\alpha_2(h)= w_2-w_3 ,\\
\alpha_3(h)= \frac12(2w_4-\imag w_5-\imag w_6),\\
\alpha_4(h)= w_3-w_4,\\\vspace{0.07cm}
\alpha_5(h)=  \frac12(2w_4+\imag w_5+\imag w_6) ,\\
\alpha_6(h)=  \frac12(w_1-w_2-w_3-w_4-\imag(2w_5-w_6)) .
\end{array}
$$
It is straightforward to check that $B$ is a basis of the root system $\Delta$, since $\Delta\subset \sum_{i=1}^6\mathbb{Z}_{\ge0}\alpha_i\cup\sum_{i=1}^6\mathbb{Z}_{\le0}\alpha_i$. 
In this occasion all the roots are noncompact, 
  $\Delta_0=\{\alpha\in \Delta\mid \alpha( \mathfrak{a})=0\}=\{\alpha\in \Delta\mid \alpha(h_i)=0\quad\forall i=1,\dots,4\}=\emptyset$, so that  $B_0=\emptyset$ and all the nodes are white.
Besides $\overline{\alpha}_1=\overline{\alpha}_6\colon \sum_{i=1}^4w_ih_i\mapsto \frac12(w_1-w_2-w_3-w_4)$ and
$\overline{\alpha}_3=\overline{\alpha}_5\colon \sum_{i=1}^4w_ih_i\mapsto  w_4$. Hence, the Satake diagram is:
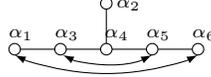
\vspace{0.2cm}

\begin{center}\begin{picture}(23,5)(4,-0.5)
\put(5,0){\circle{1}} \put(9,0){\circle{1}} \put(13,0){\circle{1}}
\put(17,0){\circle{1}} \put(21,0){\circle{1}}
\put(13,4){\circle{1}}
\put(5.5,0){\line(1,0){3}}
\put(9.5,0){\line(1,0){3}} \put(13.5,0){\line(1,0){3}} %\put(13,4.5){\line(0,1){3}}
\put(17.5,0){\line(1,0){3}}
\put(13,0.5){\line(0,1){3}}
\put(4.5,1){$\scriptstyle \alpha_1$} \put(8.5,1){$\scriptstyle \alpha_3$}
\put(12.85,1){$\scriptstyle \alpha_4$} \put(16.5,1){$\scriptstyle \alpha_5$}
\put(20.5,1){$\scriptstyle \alpha_6$} \put(13.9,3.6){$\scriptstyle \alpha_2$}
\cbezier (6,-1)(11,-2.5)(15,-2.5)(20,-1)% cubic Bezier curve
 \put(20,-1){\vector(2,1){1}}
 \put(6,-1){\vector(-2,1){1}}
 \cbezier (10.2,-1)(12.5,-1.5)(14.5,-1.5)(16.0,-1)% cubic Bezier curve
 \put(16.0,-1){\vector(2,1){1}}
 \put(10.2,-1){\vector(-2,1){1}}
 \end{picture}\end{center}\vspace{0.4cm}

Finally, the set of restricted roots $\Sigma$ is just $\Phi$, that is, a root system of type $F_4$ with basis $\{\overline{\alpha}_1,\overline{\alpha}_2,\overline{\alpha}_3,\overline{\alpha}_4 \}$. And
$\{\alpha\in \Delta\mid\overline{\alpha}=\overline{\alpha}_i\}=\{\alpha_i\}$ if $i=2,4$, so that the multiplicities are given by
$$
\begin{array}{ll}
m_{\overline{\alpha}_2 }=1, &m_{\overline{\alpha}_4 }= 1,\\
m_{\overline{\alpha}_1 }=2, &m_{\overline{\alpha}_3 }= 2.
\end{array}
$$

%%%%%%%%%%%%%%%%%%%%%%%%%%%%%%%%%%%%
%\section{Conclusions}  ¿Las ecribo o no?

%%%%%%%%%%%%%%%%%%%%%%%%%%%%%%%%%%%%%%%%%%%%%%%%%%%%%%%%%%%%%%%%%%%%%%%%%%%%%%%%%%%%%%%%%%%%%%%%%%%%%%%%%%%%%%%%%%%%%%%%%%%%%%%%%%%%%%%%%%%%%%%%%%5

%%%%%%%%%%%%%%%%%%%%%%%%%%%%%%%%%%%%%%%%%%%%%%%%%%%%%%%%%%%%%%%%%%%%%%%%%%%%%%%%%%%%%%%%%%%%%%%%%%%%%%%%%%%%%%%%%%%%%%%%%%%%%%%%%%%%%%%%%%%%%%%%%
%%%%%%%%%%%%%%%%%%%%%%%%%%%%%%%%%%%%%%%%%%%%%%%%%%%%%%%%%%%%%%%%%%%%%%%%%%%%%%%%%%%%%%%%%%%%%%%%%%%%%%%%%%%%%%%%%%%%%%%%%%%%%%%%%%%%%%%%%%%%%%%%%%%%%
\end{document}